# DISCUSSION PAPER

## ANALYSIS OF VARIANCE—WHY IT IS MORE IMPORTANT THAN EVER[1]

### By Andrew Gelman

*Columbia University*

Analysis of variance (ANOVA) is an extremely important method in exploratory and confirmatory data analysis. Unfortunately, in complex problems (e.g., split-plot designs), it is not always easy to set up an appropriate ANOVA. We propose a hierarchical analysis that automatically gives the correct ANOVA comparisons even in complex scenarios. The inferences for all means and variances are performed under a model with a separate batch of effects for each row of the ANOVA table.

We connect to classical ANOVA by working with finite-sample variance components: fixed and random effects models are characterized by inferences about existing levels of a factor and new levels, respectively. We also introduce a new graphical display showing inferences about the standard deviations of each batch of effects.

We illustrate with two examples from our applied data analysis, first illustrating the usefulness of our hierarchical computations and displays, and second showing how the ideas of ANOVA are helpful in understanding a previously fit hierarchical model.

**1. Is ANOVA obsolete?** What is the analysis of variance? Econometricians see it as an uninteresting special case of linear regression. Bayesians see it as an inflexible classical method. Theoretical statisticians have supplied many mathematical definitions [see, e.g., Speed (1987)]. Instructors see it as one of the hardest topics in classical statistics to teach, especially

Received November 2002; revised November 2003.

[1]Supported in part by the National Science Foundation with Young Investigator Award DMS-97-96129 and Grants SBR-97-08424, SES-99-87748 and SES-00-84368. A version of this paper was originally presented as a special Invited Lecture for the Institute of Mathematical Statistics.

in its more elaborate forms such as split-plot analysis. We believe, however, that the ideas of ANOVA are useful in many applications of statistics. For the purpose of this paper, we identify ANOVA with the structuring of parameters into batches—that is, with variance components models. There are more general mathematical formulations of the analysis of variance, but this is the aspect that we believe is most relevant in applied statistics, especially for regression modeling.

We shall demonstrate how many of the difficulties in understanding and computing ANOVAs can be resolved using a hierarchical Bayesian framework. Conversely, we illustrate how thinking in terms of variance components can be useful in understanding and displaying hierarchical regressions. With hierarchical (multilevel) models becoming used more and more widely, we view ANOVA as more important than ever in statistical applications.

Classical ANOVA for balanced data does three things at once:

1. As exploratory data analysis, an ANOVA is an organization of an additive data decomposition, and its sums of squares indicate the variance of each component of the decomposition (or, equivalently, each set of terms of a linear model).
2. Comparisons of mean squares, along with F-tests [or F-like tests; see, e.g., Cornfield and Tukey (1956)], allow testing of a nested sequence of models.
3. Closely related to the ANOVA is a linear model fit with coefficient estimates and standard errors.

Unfortunately, in the classical literature there is some debate on how to perform ANOVA in complicated data structures with nesting, crossing and lack of balance. In fact, given the multiple goals listed above, it is not at all obvious that a procedure recognizable as "ANOVA" should be possible at all in general settings [which is perhaps one reason that Speed (1987) restricts ANOVA to balanced designs].

In a linear regression, or more generally an additive model, ANOVA represents a batching of effects, with each row of the ANOVA table corresponding to a set of predictors. We are potentially interested in the individual coefficients and also in the variance of the coefficients in each batch. Our approach is to use variance components modeling for all rows of the table, even for those sources of variation that have commonly been regarded as fixed effects. We thus borrow many ideas from the classical variance components literature.

As we show in Section 2 of this paper, least-squares regression solves some ANOVA problems but has trouble with hierarchical structures [see also Gelman (2000)]. In Sections 3 and 4 we present a more general hierarchical regression approach that works in all ANOVA problems in which effects are structured into exchangeable batches, following the approach of Sargent and Hodges (1997). In this sense, ANOVA is indeed a special case of linear regression, but only if hierarchical models are used. In fact, the batching of effects in a hierarchical model has an exact counterpart in the rows of the analysis of variance table. Section 5 presents a new analysis of variance table that we believe more directly addresses the questions of interest in linear models, and Section 6 discusses the distinction between fixed and random effects. We present two applied examples in Section 7 and conclude with some open problems in Section 8.

**2. ANOVA and linear regression.** We begin by reviewing the benefits and limitations of classical nonhierarchical regression for ANOVA problems.

2.1. *ANOVA and classical regression*: *good news.* It is well known that many ANOVA computations can be performed using linear regression computations, with each row of the ANOVA table corresponding to the variance of a corresponding set of regression coefficients.

2.1.1. *Latin square.* For a simple example, consider a Latin square with five treatments randomized to a $5 \times 5$ array of plots. The ANOVA regression has 25 data points and the following predictors: one constant, four rows, four columns and four treatments, with only four in each batch because, if all five were included, the predictors would be collinear. (Although not necessary for understanding the mathematical structure of the model, the details of counting the predictors and checking for collinearity are important in actually implementing the regression computation and are relevant to the question of whether ANOVA can be computed simply using classical regression. As we shall discuss in Section 3.1, we ultimately will find it more helpful to include all five predictors in each batch using a hierarchical regression framework.)

For each of the three batches of variables in the Latin square problem, the variance of the $J = 5$ underlying coefficients can be estimated using the basic variance decomposition formula, where we use the notation $\mathrm{var}_{j=1}^{J}$ for the sample variance of $J$ items:

$$\mathrm{E}(\text{variance between the } \hat{\beta}_j\text{'s}) = \text{variance between the true } \beta_j\text{'s}$$

$$+ \text{estimation variance,}$$

(1)

$$\mathrm{E}(\mathrm{var}_{j=1}^{J} \hat{\beta}_j) = \mathrm{var}_{j=1}^{J} \beta_j + \mathrm{E}(\mathrm{var}(\hat{\beta}_j | \beta_j)),$$

$$\mathrm{E}(V(\hat{\beta})) = V(\beta) + V_{\text{estimation}}.$$

One can compute $V(\hat{\beta})$ and an estimate of $V_{\text{estimation}}$ directly from the coefficient estimates and standard errors, respectively, in the linear regression output, and then use the simple unbiased estimate,

$$(2) \qquad\qquad \widehat{V}(\beta) = V(\hat{\beta}) - \widehat{V}_{\text{estimation}}.$$

[More sophisticated estimates of variance components are possible; see, e.g., Searle, Casella and McCulloch (1992).] An F-test for null treatment effects corresponds to a test that $V(\beta) = 0$.

Unlike in the usual ANOVA setup, here we do not need to decide on the comparison variances (i.e., the denominators for the F-tests). The regression automatically gives standard errors for coefficient estimates that can directly be input into $\widehat{V}_{\text{estimation}}$ in (2).

2.1.2. *Comparing two treatments.* The benefits of the regression approach can be further seen in two simple examples. First, consider a simple experiment with 20 units completely randomized to two treatments, with each treatment applied to 10 units. The regression has 20 data points and two predictors: one constant and one treatment indicator (or no constant and two treatment indicators). Eighteen degrees of freedom are available to estimate the residual variance, just as in the corresponding ANOVA.

Next, consider a design with 10 pairs of units, with the two treatments randomized within each pair. The corresponding regression analysis has 20 data points and 11 predictors: one constant, one indicator for treatment and nine indicators for pairs, and, if you run the regression, the standard errors for the treatment effect estimates are automatically based on the nine degrees of freedom for the within-pair variance.

The different analyses for paired and unpaired designs are confusing for students, but here they are clearly determined by the principle of including in the regression all the information used in the design.

2.2. *ANOVA and classical regression: bad news.* Now we consider two examples where classical nonhierarchical regression *cannot* be used to automatically get the correct answer.

2.2.1. *A split-plot Latin square.* Here is the form of the analysis of variance table for a $5 \times 5 \times 2$ split-plot Latin square: a standard experimental design but one that is complicated enough that most students analyze it incorrectly unless they are told where to look it up. (We view the difficulty of teaching these principles as a sign of the awkwardness of the usual theoretical framework of these ideas rather than a fault of the students.)

In this example, there are 25 plots with five full-plot treatments (labeled A, B, C, D, E), and each plot is divided into two subplots with subplot varieties (labeled 1 and 2). As is indicated by the horizontal lines in the ANOVA table, the main-plot residual mean squares should be used for the

| Source | df |
|---|---|
| row | 4 |
| column | 4 |
| (A, B, C, D, E) | 4 |
| plot | 12 |
| (1, 2) | 1 |
| row × (1, 2) | 4 |
| column × (1, 2) | 4 |
| (A, B, C, D, E) × (1, 2) | 4 |
| plot × (1, 2) | 12 |

main-plot effects and the sub-plot residual mean squares for the sub-plot effects.

It is not hard for a student to decompose the 49 degrees of freedom to the rows in the ANOVA table; the tricky part of the analysis is to know which residuals are to be used for which comparisons.

What happens if we input the data into the `aov` function in the statistical package `S-Plus`? This program uses the linear-model fitting routine `lm`, as one might expect based on the theory that analysis of variance is a special case of linear regression. [E.g., Fox (2002) writes, "It is, from one point of view, unnecessary to consider analysis of variance models separately from the general class of linear models."] Figure 1 shows three attempts to fit the split-plot data with `aov`, only the last of which worked. We include this not to disparage `S-Plus` in any way but just to point out that ANOVA can be done in many ways in the classical linear regression framework, and not all these ways give the correct answer.

At this point, we seem to have the following "method" for analysis of variance: first, recognize the form of the problem (e.g., split-plot Latin square); second, look it up in an authoritative book such as Snedecor and Cochran (1989) or Cochran and Cox (1957); third, perform the computations, using the appropriate residual mean squares. This is unappealing for practice as well as teaching and in addition contradicts the idea that, "If you know linear regression, you know ANOVA."

2.2.2. *A simple hierarchical design.* We continue to explore the difficulties of regression for ANOVA with a simple example. Consider an experiment on four treatments for an industrial process applied to 20 machines (randomly divided into four groups of 5), with each treatment applied six times independently on each of its five machines. For simplicity, we assume no systematic time effects, so that the six measurements are simply replications. The ANOVA table is then

There are no rows for just "machine" or "measurement" because the design is fully nested.

```
> summary (aov (data ~ rows + columns + tABCDE + plots +
  t12*rows + t12*columns + t12*tABCDE + t12*plots))
            Df Sum Sq Mean Sq F value    Pr(>F)
rows         4 288.48   72.12  4.0283 0.0268475 *
columns      4 389.48   97.37  5.4387 0.0098253 **
tABCDE       4 702.28  175.57  9.8066 0.0009245 ***
plots       12 308.04   25.67  1.4338 0.2710432
t12          1 332.82  332.82 18.5898 0.0010110 **
rows:t12     4  74.08   18.52  1.0344 0.4291297
columns:t12  4  96.68   24.17  1.3500 0.3079352
tABCDE:t12   4  57.08   14.27  0.7971 0.5496092
Residuals   12 214.84   17.90

> summary (aov (data ~ rows + columns + tABCDE +
  t12*rows + t12*columns + t12*tABCDE + t12*plots + Error(plots)))
Error: plots
            Df Sum Sq Mean Sq F value Pr(>F)
rows         4 288.48   72.12  7.3592 0.2689
columns      4 389.48   97.37  9.9357 0.2331
tABCDE       4 702.28  175.57 17.9153 0.1752
plots       11 298.24   27.11  2.7666 0.4401
Residuals    1   9.80    9.80

Error: Within
            Df Sum Sq Mean Sq F value Pr(>F)
t12          1 332.82  332.82  7.3960 0.2243
rows:t12     4  74.08   18.52  0.4116 0.8059
columns:t12  4  96.68   24.17  0.5371 0.7559
tABCDE:t12   4  57.08   14.27  0.3171 0.8496
t12:plots   11 169.84   15.44  0.3431 0.8842
Residuals    1  45.00   45.00

> summary (aov (data ~ rows + columns + tABCDE +
  t12*rows + t12*columns + t12*tABCDE + Error(plots)))
Error: plots
            Df Sum Sq Mean Sq F value    Pr(>F)
rows         4 288.48   72.12  2.8095 0.073984 .
columns      4 389.48   97.37  3.7931 0.032271 *
tABCDE       4 702.28  175.57  6.8395 0.004154 **
Residuals   12 308.04   25.67

Error: Within
            Df Sum Sq Mean Sq F value    Pr(>F)
t12          1 332.82  332.82 18.5898 0.001011 **
rows:t12     4  74.08   18.52  1.0344 0.429130
columns:t12  4  96.68   24.17  1.3500 0.307935
tABCDE:t12   4  57.08   14.27  0.7971 0.549609
Residuals   12 214.84   17.90
```

FIG. 1. *Three attempts at running the* aov *command in* S-Plus. *Only the last gave the correct comparisons. This is not intended as a criticism of* S-Plus; *in general, classical ANOVA requires careful identification of variance components in order to give the correct results with hierarchical data structures.*

| Source | df |
|---|---|
| treatment | 3 |
| treatment × machine | 16 |
| treatment × machine × measurement | 100 |

Without knowing ANOVA, is it possible to get appropriate inferences for the treatment effects using linear regression? The averages for the treatments $i = 1, \ldots, 4$ can be written in two ways:

$$(3) \qquad \bar{y}_{i\cdot\cdot} = \tfrac{1}{30} \sum_{j=1}^{5} \sum_{k=1}^{6} y_{ijk}$$

or

$$(4) \qquad \bar{y}_{i\cdot\cdot} = \tfrac{1}{5} \sum_{j=1}^{5} \bar{y}_{ij\cdot}$$

Formula (3) uses all the data and suggests a standard error based on 29 degrees of freedom for each treatment, but this would ignore the nesting in the design. Formula (4) follows the design and suggests a standard error based on the four degrees of freedom from the five machines for each treatment.

Formulas (3) and (4) give the same estimated treatment effects but imply different standard errors and different ANOVA F-tests. If there is any chance of machine effects, the second analysis is standard. However, to do this you must know to base your uncertainties on the "treatment × machine" variance, not the "treatment × machine × measurement" variance. An automatic ANOVA program must be able to automatically correctly choose this comparison variance.

Can this problem be solved using least-squares regression on the 120 data points? The simplest regression uses four predictors—one constant term and three treatment indicators—with 116 residual degrees of freedom. This model gives the wrong residual variance: we want the between-machine, not the between-measurement, variance.

Since the machines are used in the design, they should be included in the analysis. This suggests a model with 24 predictors: one constant, three treatment indicators, and 20 machine indicators. But these predictors are collinear, so we must eliminate four of the machine indicators. Unfortunately, the standard errors of the treatment effects in this model are estimated using the within-machine variation, which is still wrong. The problem becomes even more difficult if the design is unbalanced.

The appropriate analysis, of course, is to include the 20 machines as a variance component, which classically could be estimated using REML

(treating the machine effects as missing data) or using regression without machine effects but with a block-structured covariance matrix with intra-class correlation estimated from data. In a Bayesian context the machine effects would be estimated with a population distribution whose variance is estimated from data, as we discuss in general in the next section. In any case, we would like to come at this answer simply by identifying the important effects—treatments and machines—without having to explicitly recognize the hierarchical nature of the design, in the same way that we would like to be able to analyze split-plot data without the potential mishaps illustrated in Figure 1.

## 3. ANOVA using hierarchical regression.

3.1. *Formulation as a regression model.* We shall work with linear models, with the "analysis of variance" corresponding to the batching of effects into "sources of variation," and each batch corresponding to one row of the ANOVA table. This is the model of Sargent and Hodges (1997). We use the notation $m = 1, \ldots, M$ for the rows of the table. Each row $m$ represents a batch of $J_m$ regression coefficients $\beta_j^{(m)}$, $j = 1, \ldots, J_m$. We denote the $m$th subvector of coefficients as $\beta^{(m)} = (\beta_1^{(m)}, \ldots, \beta_{J_m}^{(m)})$ and the corresponding classical least-squares estimate as $\hat{\beta}^{(m)}$. These estimates are subject to $c_m$ linear constraints, yielding $(df)_m = J_m - c_m$ degrees of freedom. We label the constraint matrix as $C^{(m)}$, so that $C^{(m)} \hat{\beta}^{(m)} = 0$ for all $m$. For notational convenience, we label the grand mean as $\beta_1^{(0)}$, corresponding to the (invisible) zeroth row of the ANOVA table and estimated with no linear constraints.

The linear model is fit to the data points $y_i$, $i = 1, \ldots, n$, and can be written as

$$(5) \qquad y_i = \sum_{m=0}^{M} \beta_{j_i^m}^{(m)},$$

where $j_i^m$ indexes the appropriate coefficient $j$ in batch $m$ corresponding to data point $i$. Thus, each data point pulls one coefficient from each row in the ANOVA table. Equation (5) could also be expressed as a linear regression model with a design matrix composed entirely of 0's and 1's. The coefficients $\beta_j^M$ of the last row of the table correspond to the residuals or error term of the model. ANOVA can also be applied more generally to regression models (or to generalized linear models), in which case we could have any design matrix $X$, and (5) would be generalized to

$$(6) \qquad y_i = \sum_{m=0}^{M} \sum_{j=1}^{J_m} x_{ij}^{(m)} \beta_j^{(m)}.$$

The essence of analysis of variance is in the structuring of the coefficients into batches—hence the notation $\beta_j^{(m)}$—going beyond the usual linear model formulation that has a single indexing of coefficients $\beta_j$. We assume that the structure (5), or the more general regression parameterization (6), has already been constructed using knowledge of the data structure. To use ANOVA terminology, we assume the sources of variation have already been set, and our goal is to perform inference for each variance component.

We shall use a hierarchical formulation in which each batch of regression coefficients is modeled as a sample from a normal distribution with mean 0 and its own variance $\sigma_m^2$:

$$(7) \quad \beta_j^{(m)} \sim \mathrm{N}(0, \sigma_m^2) \qquad \text{for } j = 1, \ldots, J_m \text{ for each batch } m = 1, \ldots, M.$$

We follow the notation of Nelder (1977, 1994) by modeling the underlying $\beta$ coefficients as unconstrained, unlike the least-squares estimates. Setting the variances $\sigma_m^2$ to $\infty$ and constraining the $\beta_j^{(m)}$'s yields classical least-squares estimates.

Model (7) corresponds to exchangeability of each set of factor levels, which is a form of partial exchangeability or invariance of the entire set of cell means [see Aldous (1981)]. We do not mean to suggest that this model is universally appropriate for data but rather that it is often used, explicitly or implicitly, as a starting point for assessing the relative importance of the effects $\beta$ in linear models structured as in (5) and (6). We discuss nonexchangeable models in Section 8.3.

One measure of the importance of each row or "source" in the ANOVA table is the standard deviation of its constrained regression coefficients, which we denote

$$(8) \qquad s_m = \sqrt{\frac{1}{(df)_m} \beta^{(m)T} [I - C^{(m)} (C^{(m)T} C^{(m)})^{-1} C^{(m)T}] \beta^{(m)}},$$

where $\beta^{(m)}$ is the vector of coefficients in batch $m$ and $C^{(m)}$ is the $c_m \times J_m$ full rank matrix of constraints (for which $C^{(m)} \beta^{(m)} = 0$). Expression (8) is just the mean square of the coefficients' residuals after projection to the constraint space. We divide by $(df)_m = J_m - c_m$ rather than $J_m - 1$ because multiplying by $C^{(m)}$ induces $c_m$ linear constraints.

Variance estimation is often presented in terms of the superpopulation standard deviations $\sigma_m$, but in our ANOVA summaries we focus on the finite-population quantities $s_m$ for reasons discussed in Section 3.5. However, for computational reasons the parameters $\sigma_m$ are useful intermediate quantities to estimate.

3.2. *Batching of regression coefficients.*   Our general solution to the ANOVA problem is simple: we treat *every* row in the table as a batch of "random effects"; that is, a set of regression coefficients drawn from a distribution with mean 0 and some standard deviation to be estimated from the data. The mean of 0 comes naturally from the ANOVA decomposition structure (pulling out the grand mean, main effects, interactions and so forth), and the standard deviations are simply the magnitudes of the variance components corresponding to each row of the table. For example, we can write the simple hierarchical design of Section 2.2.2 as

| Source | Number of coefficients | Standard deviation |
|---|---|---|
| treatment | 4 | $s_1$ |
| treatment × machine | 20 | $s_2$ |
| treatment × machine × measurement | 120 | $s_3$ |

Except for our focus on $s$ rather than $\sigma$, this is the approach recommended by Box and Tiao (1973) although computational difficulties made it difficult to implement at that time.

The primary goal of ANOVA is to estimate the variance components (in this case, $s_1, s_2, s_3$) and compare them to zero and to each other. The secondary goal is to estimate (and summarize the uncertainties in) the individual coefficients, especially, in this example, the four treatment effects. From the hierarchical model the coefficient estimates will be pulled toward zero, with the amount of shrinkage determined by the estimated variance components. But, more importantly, the variance components and standard errors are estimated from the data, without any need to specify comparisons based on the design. Thus, the struggles of Section 2.2 are avoided, and (hierarchical) linear regression can indeed be used to compute ANOVA automatically, once the rows of the table (the sources of variation) have been specified.

For another example, the split-plot Latin square looks like

This is automatic, based on the principle that all variables in the design be included in the analysis. Setting up the model in this way, with all nine variance components estimated, automatically gives the correct comparisons (e.g., uncertainties for comparisons between treatments A, B, C, D, E will be estimated based on main-plot variation and uncertainties for varieties 1, 2 will be estimated based on sub-plot variation).

| Source | Number of coefficients | Standard deviation |
|---|---|---|
| row | 5 | $s_1$ |
| column | 5 | $s_2$ |
| (A, B, C, D, E) | 5 | $s_3$ |
| plot | 25 | $s_4$ |
| (1, 2) | 2 | $s_5$ |
| row × (1, 2) | 10 | $s_6$ |
| column × (1, 2) | 10 | $s_7$ |
| (A, B, C, D, E) × (1, 2) | 10 | $s_8$ |
| plot × (1, 2) | 50 | $s_9$ |

3.3. *Getting something for nothing*? At this point we seem to have a paradox. In classical ANOVA, you (sometimes) need to know the design in order to select the correct analysis, as in the examples in Section 2.2. But the hierarchical analysis does it automatically. How can this be? How can the analysis "know" how to do the split-plot analysis, for example, without being "told" that the data come from a split-plot design?

The answer is in two parts. First, as with the classical analyses, we require that the rows of the ANOVA be specified by the modeler. In the notation of (5) and (6), the user must specify the structuring or batching of the linear parameters $\beta$. In the classical analysis, however, this is not enough, as discussed in Section 2.2.

The second part of making the hierarchical ANOVA work is that the information from the design is encoded in the design matrix of the linear regression [as shown by Nelder (1965a, b) and implemented in the software Genstat]. For example, the nesting in the example of Section 2.2.2 is reflected in the collinearity of the machine indicators within each treatment. The automatic encoding is particularly useful in incomplete designs where there is no simple classical analysis.

From a linear-modeling perspective, classical nonhierarchical regression has a serious limitation: each batch of parameters (corresponding to each row of the ANOVA table) must be included with no shrinkage (i.e., $\sigma_m = \infty$) or excluded ($\sigma_m = 0$), with the exception of the last row of the table, whose variance can be estimated. In the example of Section 2.2.2, we must either include the machine effects unshrunken or ignore them, and neither approach gives the correct analysis. The hierarchical model works automatically because it allows finite nonzero values for all the variance components.

The hierarchical regression analysis is based on the model of exchangeable effects within batches, as expressed in model (7), which is not necessarily the best analysis in any particular application. For example, Besag and Higdon (1999) recommend using spatial models (rather than exchangeable row and column effects) for data such as in the split-plot experiment described previously. Here we are simply trying to understand why, when given the standard assumptions underlying the classical ANOVA, the hierarchical analysis automatically gives the appropriate inferences for the variance components without the need for additional effort of identifying appropriate error terms for each row of the table.

3.4. *Classical and Bayesian interpretations.* We are most comfortable interpreting the linear model in a Bayesian manner, that is, with a joint probability distribution on all unknown parameters. However, our recommended hierarchical approach can also be considered classically, in which case the regression coefficients are considered as random variables (and thus are "predicted") and the variance components are considered as parameters (and thus "estimated"); see Robinson (1991) and Gelman, Carlin, Stern and Rubin [(1995), page 380]. The main difference between classical and Bayesian methods here is between using a point estimate for the variance parameters or including uncertainty distributions. Conditional on the parameters $\sigma_m$, the classical and Bayesian inferences for the linear parameters $\beta_j^m$ are identical in our ANOVA models. In either case, the individual regression coefficients are estimated by linear unbiased predictors or, equivalently, posterior means, balancing the direct information on each parameter with the shrinkage from the batch of effects. There will be more shrinkage for batches of effects whose standard deviations $\sigma_m$ are near zero, which will occur for factors that contribute little variation to the data.

When will it make a practical difference to estimate variance parameters Bayesianly rather than with point estimates? Only when these variances are hard to distinguish from 0. For example, Figure 2 shows the posterior distribution of the hierarchical standard deviation from an example of Rubin (1981) and Gelman, Carlin, Stern and Rubin [(1995), Chapter 5]. The data are consistent with a standard deviation of 0, but it could also be as high as 10 or 20. Setting the variance parameter to zero in such a situation is generally *not* desirable because it would lead to falsely precise estimates of the $\beta_j^{(m)}$'s. Setting the variance to some nonzero value would require additional work which, in practice, would not be done since it would offer no advantages over Bayesian posterior averaging.

It might be argued that such examples—in which the maximum likelihood estimate of the hierarchical variance is at or near zero—are pathological and unlikely to occur in practice. But we would argue that such situations will be common in ANOVA settings, for two reasons. First, when studying the many rows of a large ANOVA table, we expect (in fact, we hope) to see various near-zero variances at higher levels of interaction. After all, one of the pur-poses of an ANOVA decomposition is to identify the important main effects and interactions in a complex data set [see Sargent and Hodges (1997)]. Nonsignificant rows of the ANOVA table correspond to variance components that are statistically indistinguishable from zero. Our second reason for expecting to see near-zero variance components is that, as informative covariates are added to a linear model, hierarchical variances decrease until it is no longer possible to add more information [see Gelman (1996)].

When variance parameters are not well summarized by point estimates, Bayesian inferences are sensitive to the prior distribution. For our basic ANOVA computations we use noninformative prior distributions of the form $p(\sigma_m) \propto 1$ (which can be considered as a degenerate case of the inverse-gamma family, as we discuss in Section 4.2). We further discuss the issue of near-zero variance components in Section 8.2.

3.5. *Superpopulation and finite-population variances.* For each row $m$ of an ANOVA table, there are two natural variance parameters to estimate: the *superpopulation* standard deviation $\sigma_m$ and the *finite-population* standard deviation $s_m$ as defined in (8). The superpopulation standard deviation characterizes the uncertainty for predicting a new coefficient from batch $m$, whereas the finite-population standard deviation describes the existing $J_m$ coefficients. The two variances can be given the same point estimate—in classical unbiased estimation $\mathrm{E}(s_m^2|\sigma_m^2) = \sigma_m^2$, and in Bayesian inference with a noninformative prior distribution (see Section 4.2) the conditional posterior mode of $\sigma_m^2$ given all other parameters in the model is $s^2$. The superpopulation variance has more uncertainty, however.

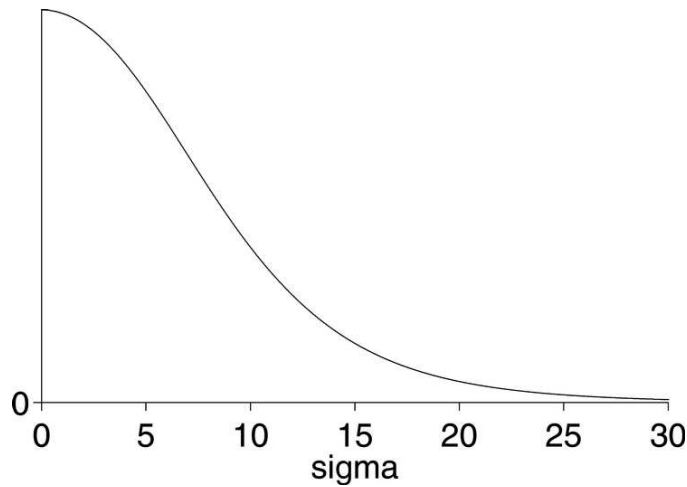

Fig. 2. *Illustration of the difficulties of point estimation for variance components. Pictured is the marginal posterior distribution for a hierarchical standard deviation parameter from Rubin (1981) and Gelman, Carlin, Stern and Rubin [(1995), Chapter 5]. The simplest point estimate, the posterior mode or REML estimate, is zero, but this estimate is on the extreme of parameter space and would cause the inferences to understate the uncertainties in this batch of regression coefficients.*

To see the difference between the two variances, consider the extreme case in which $J_m = 2$ [and so $(df)_m = 1$] and a large amount of data is available in both groups. Then the two parameters $\beta_1^{(m)}$ and $\beta_2^{(m)}$ will be estimated accurately and so will $s_m^2 = (\beta_1^{(m)} - \beta_2^{(m)})^2/2$. The superpopulation variance $\sigma_m^2$, on the other hand, is only being estimated by a measurement that is proportional to a $\chi^2$ with one degree of freedom. We know much about the two parameters $\beta_1^{(m)}, \beta_2^{(m)}$ but can say little about others from their batch.

As we discuss in Section 6, we believe that much of the literature on fixed and random effects can be fruitfully reexpressed in terms of finite-population and superpopulation inferences. In some contexts (e.g., obtaining inference for the 50 U.S. states) the finite population seems more meaningful, whereas in others (e.g., subject-level effects in a psychological experiment) interest clearly lies in the superpopulation.

To keep connection with classical ANOVA, which focuses on a description— a variance decomposition—of an existing dataset, we focus on finite-population variances $s_m^2$. However, as an intermediate step in any computation—classical or Bayesian—we perform inferences about the superpopulation variances $\sigma_m^2$.

## 4. Inference for the variance components.

4.1. *Classical inference.* Although we have argued that hierarchical models are best analyzed using Bayesian methods, we discuss classical computations first, partly because of their simplicity and partly to connect to the vast literature on the estimation of variance components [see, e.g., Searle, Casella and McCulloch (1992)]. The basic tool is the method of moments. We can first estimate the superpopulation variances $\sigma_m^2$ and their approximate uncertainty intervals, then go back and estimate uncertainty intervals for the finite-population variances $s_m^2$. Here we are working with the additive model (5) rather than the general regression formulation (6).

The estimates for the parameters $\sigma_m^2$ are standard and can be expressed in terms of classical ANOVA quantities, as follows. The sum of squares for row $m$ is the sum of the squared coefficient estimates corresponding to the $n$ data points,

$$SS_m = \sum_{i=1}^n (\hat{\beta}_{j_i^m}^{(m)})^2,$$

and can also be written as a weighted sum of the squared coefficient estimates for that row,

$$SS_m = n \sum_{j=1}^{J_m} w_j (\hat{\beta}_j^{(m)})^2,$$

where the weights $w_j$ sum to 1, and

for balanced designs:   $SS_m = \dfrac{n}{J_m} \sum_{j=1}^{J_m} (\hat{\beta}_j^{(m)})^2.$

The mean square is the sum of squares divided by degrees of freedom,

$$MS_m = SS_m / (df)_m$$

and

for balanced designs:   $MS_m = \dfrac{n}{J_m (df)_m} \sum_{j=1}^{J_m} (\hat{\beta}_j^{(m)})^2.$

The all-important expected mean square, $EMS_m$, is the expected contribution of sampling variance to $MS_m$, and it is also $\mathrm{E}(MS_m)$ under the null hypothesis that the coefficients $\beta_j^{(m)}$ are all equal to zero. Much of the classical literature is devoted to determining $EMS_m$ under different designs and different assumptions, and computing or approximating the F-ratio, $MS_m / EMS_m$, to assess statistical significance.

We shall proceed in a slightly different direction. First, we compute $EMS_m$ under the general model allowing all other variance components in the model to be nonzero. (This means that, in general, $EMS_m$ depends on variance components estimated lower down in the ANOVA table.) Second, we use the expected mean square as a tool to *estimate* variance components, not to test their statistical significance. Both these steps follow classical practice for random effects; our only innovation is to indiscriminately apply them to *all* the variance components in a model, and to follow this computation with an estimate of the uncertainty in the finite-population variances $s_m^2$.

We find it more convenient to work with not the sums of squares or mean squares but with the variances of the batches of estimated regression coefficients, which we label as

(9) $$V_m = \frac{1}{(df)_m} \sum_{j=1}^{J_m} (\hat{\beta}_j^{(m)})^2.$$

$V_m$ can be considered a variance since for each row the $J_m$ effect estimates $\hat{\beta}^{(m)}$ have several linear constraints [with $(df)_m$ remaining degrees of freedom] and must sum to 0. [For the "zeroth" row of the table, we define $V_0 = (\hat{\beta}_1^{(0)})^2$, the square of the estimated grand mean in the model.] For each row of the table,

$$\text{for balanced designs:} \quad V_m = \frac{J_m}{n} MS_m.$$

We start by estimating the superpopulation variances $\sigma_m^2$, and the constrained method-of-moments estimator is based on the variance-decomposition identity [see (1)]

$$\mathrm{E}(V_m) = \sigma_m^2 + EV_m,$$

where $EV_m$ is the contribution of sampling variance to $V_m$, that is, the expected value of $V_m$ if $\sigma_m$ were equal to 0. $EV_m$ in turn depends on other variance components in the model, and

$$\text{for balanced designs:} \quad EV_m = \frac{J_m}{n} EMS_m.$$

The natural estimate of the underlying variance is then

$$(10) \qquad\qquad \hat{\sigma}_m^2 = \max(0, V_m - \widehat{EV}_m).$$

The expected value $\widehat{EV}_m$ is itself estimated based on the other variance components in the model, as we discuss shortly.

Thus, the classical hierarchical ANOVA computations reduce to estimating the expected mean squares $EMS_m$ (and thus $EV_m$) in terms of the estimated variance components $\sigma_m$. For nonbalanced designs, this can be complicated compared to the Bayesian computation as described in Section 4.2.

For balanced designs, however, simple formulas exist. We do not go through all the literature here [see, e.g., Cornfield and Tukey (1956), Green and Tukey (1960) and Plackett (1960)]. A summary is given in Searle, Casella and McCulloch [(1992), Section 4.2]. The basic idea is that, in a balanced design, the effect estimates $\hat{\beta}_j^{(m)}$ in a batch $m$ are simply averages of data, adjusted to fit a set of linear constraints. The sampling variance $\widehat{EV}_m$ in (10) can be written in terms of variances $\sigma_k^2$ for all batches $k$ representing interactions that include $m$ in the ANOVA table. We write this as

$$(11) \qquad\qquad \widehat{EV}_m = \sum_{k \in I(m)} \frac{J_m}{J_k} \sigma_k^2,$$

where $I(m)$ represents the set of all rows in the ANOVA table representing interactions that include the variables $m$ as a subset. For example, in the example in Section 2.2.2, consider the treatment effects (i.e., $m = 1$ in the ANOVA table). Here, $J_1 = 4$, $n = 120$ and $\widehat{EV}_1 = \frac{4}{20}\sigma_2^2 + \frac{4}{120}\sigma_3^2$. For another example, in the split-plot latin square in Section 2.2.1, the main-plot treatment effects are the third row of the ANOVA table ($m = 3$), and $\widehat{EV}_3 = \frac{5}{25}\sigma_4^2 + \frac{5}{10}\sigma_8^2 + \frac{5}{50}\sigma_9^2$.

For balanced designs, then, variance components can be estimated by starting at the bottom of the table (with the highest-level interaction, or residuals) and then working upwards, at each step using the appropriate variance components from lower in the table in formulas (10) and (11). In this way the variance components $\sigma_m^2$ can be estimated noniteratively. Alternatively, we can compute the moments estimator of the entire vector $\sigma^2 = (\sigma_1^2, \ldots, \sigma_M^2)$ at once by solving the linear system $V = A\hat{\sigma}^2$, where $V$ is the vector of raw row variances $V_m$ and $A$ is the square matrix with $A_{km} = \frac{J_m}{J_k}$ if $k \in I(m)$ and 0 otherwise.

The next step is to determine uncertainties for the estimated variance components. Once again, there is an extensive literature on this; the basic method is to express each estimate $\hat{\sigma}_m^2$ as a sum and difference of independent random variables whose distributions are proportional to $\chi^2$, and then to compute the variance of the estimate. The difficulty of this standard approach is in working with this combination-of-$\chi^2$ distribution.

Instead, we evaluate the uncertainties of the estimated variance components by simulation, performing the following steps 1000 times: (1) simulate uncertainty in each raw row variance $V_m$ by multiplying by a random variable of the form $(df)_m/\chi^2_{(df)_m}$, (2) solve for $\hat{\sigma}^2$ in $V = A\hat{\sigma}^2$, (3) constrain the solution to be nonnegative, and (4) compute the 50% and 95% intervals from the constrained simulation draws. This simulation has a parametric bootstrap or Bayesian flavor and is motivated by the approximate equivalence between repeated-sampling and Bayesian inferences [see, e.g., DeGroot (1970) and Efron and Tibshirani (1993)].

Conditional on the simulation for $\sigma$, we can now estimate the finite-population standard deviations $s_m$. As discussed in Section 3.5, the data provide additional information about these, and so our intervals for $s_m$ will be narrower than for $\sigma_m$, especially for variance components with few degrees of freedom. Given $\sigma$, the parameters $\beta_j^{(m)}$ have a multivariate normal distribution (in Bayesian terms, a conditional posterior distribution; in classical terms, a predictive distribution). The resulting inference for each $s_m$ can be derived from (8), computing either by simulation of the $\beta$'s or by approximation with the $\chi^2$ distribution. Finally, averaging over the simulations of $\sigma$ yields predictive inferences about the $s_m$'s.

4.2. *Bayesian inference.* To estimate the variance components using Bayesian methods, one needs a probability model for the regression coefficients $\beta_j^{(m)}$ and the variance parameters $\sigma_m$. The standard model for $\beta$'s is independent normal, as given by (7). In our ANOVA formulation (5) or (6), the regression error terms are just the highest-level interactions, $\beta_j^{(M)}$, and so the distributions (7) include the likelihood as well as the prior distribution. For generalized linear models, the likelihood can be written separately (see Section 7.2 for an example).

The conditionally conjugate hyperprior distributions for the variances can be written as scaled inverse-$\chi^2$:

$$\sigma_m^2 \sim \text{Inv-}\chi^2(\nu_m, \sigma_{0m}^2).$$

A standard noninformative prior distribution is uniform on $\sigma$, which corresponds to each $\nu_m = -1$ and $\sigma_{0m} = 0$ [see, e.g., Gelman, Carlin, Stern and Rubin (1995)]. For values of $m$ in which $J_m$ is large (i.e., rows of the ANOVA table corresponding to many linear predictors), $\sigma_m$ is essentially estimated from data. When $J_m$ is small, the flat prior distribution implies that $\sigma$ is allowed the possibility of taking on large values, which minimizes the amount of shrinkage in the effect estimates.

More generally, it would make sense to model the variance parameters $\sigma_m$ themselves, especially for complicated models with many variance components (i.e., many rows of the ANOVA table). Such models are a potential subject of future research; see Section 8.2.

With the model as set up above, the posterior distribution for the parameters $(\beta, \sigma)$ can be simulated using the Gibbs sampler, alternately updating the vector $\beta$ given $\sigma$ with linear regression, and updating the vector $\sigma$ from the independent inverse-$\chi^2$ conditional posterior distributions given $\beta$. The only trouble with this Gibbs sampler is that it can get stuck with variance components $\sigma_m$ near zero. A more efficient updating reparameterizes into vectors $\gamma$, $\alpha$ and $\tau$, which are defined as follows:

$$(12) \qquad \begin{aligned} \beta_j^{(m)} &= \alpha_m \gamma_j^{(m)}, \\ \sigma_m &= \alpha_m \tau_m. \end{aligned}$$

The model can be then expressed as

$$y = X(\alpha\gamma),$$

$$\gamma_j^{(m)} \sim \text{N}(0, \tau_m^2) \qquad \text{for each } m,$$

$$\tau_m^2 \sim \text{Inv-}\chi^2(\nu_m, \sigma_{0m}^2).$$

The auxiliary parameters $\alpha$ are given a uniform prior distribution, and then this reduces to the original model [see Boscardin (1996), Meng and van Dyk (1997), Liu, Rubin and Wu (1998), Liu and Wu (1999) and Gelman (2004)]. The Gibbs sampler then proceeds by updating $\gamma$ (using linear regression with $n$ data points and $\sum_{m=0}^{M} J_m$ predictors), $\alpha$ (linear regression with $n$ data points and $M$ predictors) and $\tau^2$ (independent inverse-$\chi^2$ distributions). The parameters in the original parameterization, $\beta$ and $\sigma$, can then be recomputed from (12) and stored at each step.

Starting points for the Bayesian computation can be adapted from the classical point estimates for $\sigma^2$ and their uncertainties from Section 4.1. The only difficulty is that the variance parameters cannot be set to exactly zero. One reasonable approach is to replace any $\sigma_m^2$ of zero by a random value between zero and $|V_m - \widehat{EV}_m|$, treating this absolute value as a rough measure of the noise level in the estimate. Generalized linear models can be computed using this Gibbs sampler with Metropolis jumping for the nonconjugate conditional densities [see, e.g., Gelman, Carlin, Stern and Rubin (1995)] or data augmentation [see Albert and Chib (1993) and Liu (2002)]. In either case, once the simulations have approximately converged and posterior simulations are available, one can construct simulation-based intervals for all the parameters and for derived quantities of interest such as the finite-population standard deviations $s_m$ defined in (8).

When we use the uniform prior density for the parameters $\sigma_m$, the posterior distributions are proper for batches $m$ with at least two degrees of freedom. However, for effects that are unique or in pairs [i.e., batches for which $(df)_m = 1$], the posterior density for the corresponding $\sigma_m$ is improper, with infinite mass in the limit $\sigma_j \to \infty$ [Gelman, Carlin, Stern and Rubin (1995), Exercise 5.8], and so the coefficients $\beta_j^{(m)}$ in these batches are essentially being estimated via maximum likelihood. This relates to the classical result that shrinkage estimation dominates least squares when estimating three or more parameters in a normal model [James and Stein (1961)].

## 5. A new ANOVA table.
There is room for improvement in the standard analysis of variance table: it is read in order to assess the relative importance of different sources of variation, but the numbers in the table do not directly address this issue. The sums of squares are a decomposition of the total sum of squares, but the lines in the table with higher sums of squares are not necessarily those with higher estimated underlying variance components. The mean square for each row has the property that, if the corresponding effects are all zero, its expectation equals that of the error mean square. Unfortunately, if these other effects are *not* zero, the mean square has no direct interpretation in terms of the model parameters. The mean square is the variance explained per parameter, which is not directly comparable to the parameters $s_m^2$ and $\sigma_m^2$, which represent underlying variance components.

Similarly, statistical significance (or lack thereof) of the mean squares is relevant; however, rows with higher F-ratios or more extreme $p$-values do *not* necessarily correspond to batches of effects with higher estimated magnitudes. In summary, the standard ANOVA table gives all sorts of information, but nothing to directly compare the listed sources of variation.

Our alternative ANOVA table presents, for each source of variation $m$, the estimates and uncertainties for $s_m$, the standard deviation of the coefficients corresponding to that row of the table. In addition to focusing on estimation rather than testing, we display the estimates and uncertainties graphically. Since the essence of ANOVA is comparing the importance of different rows of the table, it is helpful to allow direct graphical comparison, as with tabular displays in general [see Gelman, Pasarica and Dodhia (2002)]. In addition, using careful formatting, we can display this in no more space than is required by the classical ANOVA table.

Figure 3 shows an example with the split-plot data that we considered earlier. For each source of variation, the method-of-moments estimate of $s_m$ is shown by a point, with the thick and thin lines showing 50% and 95% intervals from the simulations. The point estimates are not always at the center of the intervals because of edge effects caused by the restriction that all the variance components be nonnegative. In an applied context it might make sense to use as point estimates the medians of the simulations. We display the moments estimates here to show the effects of the constrained inference in an example where uncertainty is large.

In our ANOVA table, the inferences for all the variance components are simultaneous, in contrast to the classical approach in which each variance component is tested under the model that all others, except for the error term, are zero. Thus, the two tables answer different inferential questions. We would argue that the simultaneous inference is more relevant in applications. However, if the classical $p$-values are of interest, they could be incorporated into our graphical display.

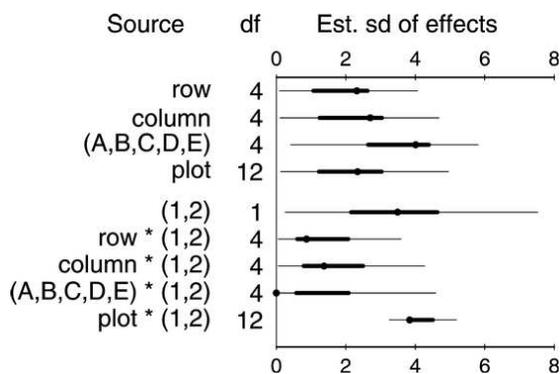

FIG. 3. ANOVA display for a split-plot latin square experiment (cf. to the classical ANOVA, which is the final table in Figure 1). The points indicate classical variance component estimates, and the bars display 50% and 95% intervals for the finite-population standard deviations $\sigma_m$. The confidence intervals are based on simulations assuming the variance parameters are nonnegative; as a result, they can differ from the point estimates, which are based on the method of moments, truncating negative estimates to zero.

**6. Fixed and random effects.**   A persistent point of conflict in the ANOVA literature is the appropriate use of fixed or random effects, an issue which we must address since we advocate treating *all* batches of effects as sets of random variables. Eisenhart ([1947](#)) distinguishes between fixed and random effects in estimating variance components, and this approach is standard in current textbooks [e.g., Kirk ([1995](#))]. However, there has been a stream of dissenters over the years; for example, Yates ([1967](#)):

> ... whether the factor levels are a random selection from some defined set (as might be the case with, say, varieties), or are deliberately chosen by the experimenter, does not affect the logical basis of the formal analysis of variance or the derivation of variance components.

Before discussing the technical issues, we briefly review what is meant by fixed and random effects. It turns out that different—in fact, incompatible—definitions are used in different contexts. [See also Kreft and de Leeuw ([1998](#)), Section 1.3.3, for a discussion of the multiplicity of definitions of fixed and random effects and coefficients, and Robinson ([1998](#)) for a historical overview.] Here we outline five definitions that we have seen:

1. Fixed effects are constant across individuals, and random effects vary. For example, in a growth study, a model with random intercepts $\alpha_i$ and fixed slope $\beta$ corresponds to parallel lines for different individuals $i$, or the model $y_{it} = \alpha_i + \beta t$. Kreft and de Leeuw [([1998](#)), page 12] thus distinguish between fixed and random coefficients.
2. Effects are fixed if they are interesting in themselves or random if there is interest in the underlying population. Searle, Casella and McCulloch [([1992](#)), Section 1.4] explore this distinction in depth.
3. "When a sample exhausts the population, the corresponding variable is *fixed*; when the sample is a small (i.e., negligible) part of the population the corresponding variable is *random*" [Green and Tukey ([1960](#))].
4. "If an effect is assumed to be a realized value of a random variable, it is called a random effect" [LaMotte ([1983](#))].
5. Fixed effects are estimated using least squares (or, more generally, maximum likelihood) and random effects are estimated with shrinkage ["linear unbiased prediction" in the terminology of Robinson ([1991](#))]. This definition is standard in the multilevel modeling literature [see, e.g., Snijders and Bosker ([1999](#)), Section 4.2] and in econometrics.

   In the Bayesian framework, this definition implies that fixed effects $\beta_j^{(m)}$ are estimated conditional on $\sigma_m = \infty$ and random effects $\beta_j^{(m)}$ are estimated conditional on $\sigma_m$ from the posterior distribution.

Of these definitions, the first clearly stands apart, but the other four definitions differ also. Under the second definition, an effect can change from fixed to random with a change in the goals of inference, even if the data and design are unchanged. The third definition differs from the others in defining a finite population (while leaving open the question of what to do with a large but not exhaustive sample), while the fourth definition makes no reference to an actual (rather than mathematical) population at all. The second definition allows fixed effects to come from a distribution, as long as that distribution is not of interest, whereas the fourth and fifth do not use any distribution for inference about fixed effects. The fifth definition has the virtue of mathematical precision but leaves unclear when a given set of effects should be considered fixed or random. In summary, it is easily possible for a factor to be "fixed" according to some of the definitions above and "random" for others. Because of these conflicting definitions, it is no surprise that "clear answers to the question 'fixed or random?' are not necessarily the norm" [Searle, Casella and McCulloch (1992), page 15].

One way to focus a discussion of fixed and random effects is to ask how inferences change when a set of effects is changed from fixed to random, with no change in the data. For example, suppose a factor has four degrees of freedom corresponding to five different medical treatments, and these are the only existing treatments and are thus considered "fixed" (according to definitions 2 and 3 above). Suppose it is then discovered that these are part of a larger family of many possible treatments, and so it is desired to model them as "random." In the framework of this paper, the inference about these five parameters $\beta_j^{(m)}$ and their finite-population and superpopulation standard deviations, $s_m$ and $\sigma_m$, will not change with the news that they actually are viewed as a random sample from a distribution of possible treatment effects. But the superpopulation variance now has an important new role in characterizing this distribution. The difference between fixed and random effects is thus not a difference in inference or computation but in the ways that these inferences will be used. Thus, we strongly *disagree* with the claim of Montgomery [(1986), page 45] that in the random effects model, "knowledge about particular [regression coefficients] is relatively useless."

We prefer to sidestep the overloaded terms "fixed" and "random" with a cleaner distinction by simply renaming the terms in definition 1 above. We define effects (or coefficients) in a multilevel model as *constant* if they are identical for all groups in a population and *varying* if they are allowed to differ from group to group. For example, the model $y_{ij} = \alpha_j + \beta x_{ij}$ (of units $i$ in groups $j$) has a constant slope and varying intercepts, and $y_{ij} = \alpha_j + \beta_j x_{ij}$ has varying slopes and intercepts. In this terminology (which we would apply at any level of the hierarchy in a multilevel model), varying effects occur in batches, whether or not the effects are interesting in themselves (definition 2), and whether or not they are a sample from a larger set (definition 3). Definitions 4 and 5 do not arise for us since we estimate all batches of effects hierarchically, with the variance components $\sigma_m$ estimated from data.

**7. Examples.** We give two examples from our own consulting and research where ANOVA has been helpful in understanding the structure of variation in a dataset. Section 7.1 describes a multilevel linear model for a full-factorial dataset, and Section 7.2 describes a multilevel logistic regression.

From a classical perspective of inference for variance components, these cases can be considered as examples of the effectiveness of automatically setting up hierarchical models with random effects for each row in the ANOVA table. From a Bayesian perspective, these examples demonstrate how the ANOVA idea—batching effects into rows and considering the importance of each batch—applies outside of the familiar context of hypothesis testing.

7.1. *A five-way factorial structure*: *Web connect times*. Data were collected by an Internet infrastructure provider on connect times—the time required for a signal to reach a specified destination—as processed by each of two different companies. Messages were sent every hour for 25 consecutive hours, from each of 45 locations to four different destinations, and the study was repeated one week later. It was desired to quickly summarize these data to learn about the importance of different sources of variation in connect times.

Figure 4 shows a classical ANOVA of logarithms of connect times using the standard factorial decomposition on the five factors: destination ("to"), source ("from"), service provider ("company"), time of day ("hour") and week. The data have a full factorial structure with no replication, so the full five-way interaction, at the bottom of the table, represents the "error" or lowest-level variability. The ANOVA reveals that all the main effects and almost all the interactions are statistically significant. However, as discussed in Section 5, it is difficult to use these significance levels, or the associated sums of squares, mean squares or F-statistics, to *compare* the importance of the different factors.

Figure 5 shows the full multilevel ANOVA display for these data. Each row shows the estimated finite-population standard deviation of the corresponding group of parameters, along with 50% and 95% uncertainty intervals. We can now immediately see that the lowest-level variation is more important in variance than any of the factors except for the main effect of the destination. `Company` has a large effect on its own and, perhaps more interestingly, in interaction with `to`, `from`, and in the three-way interaction.

The information in the multilevel display in Figure 5 is *not* simply contained in the mean squares of the classical ANOVA table in Figure 4. For example, the effects of `from * hour` have a relatively high estimated standard deviation but a relatively low mean square (see, e.g., `to * week`).

Figure 5 does not represent the end of any statistical analysis; for example, in this problem the analysis has ignored any geographical structure in the

| Source | Df | Ss | Ms | Fstat | Pvalue |
|---|---|---|---|---|---|
| to | 3 | 31193.62 | 10397.87 | 26660.68 | 0.00 |
| from | 44 | 5635.24 | 128.07 | 328.39 | 0.00 |
| company | 1 | 1027.44 | 1027.44 | 2634.40 | 0.00 |
| hour | 24 | 128.74 | 5.36 | 13.75 | 0.00 |
| week | 1 | 3.76 | 3.76 | 9.64 | 0.00 |
| | | | | | |
| to * from | 132 | 669.56 | 5.07 | 13.01 | 0.00 |
| to * company | 3 | 497.03 | 165.68 | 424.80 | 0.00 |
| to * hour | 72 | 44.00 | 0.61 | 1.57 | 0.00 |
| to * week | 3 | 14.59 | 4.86 | 12.47 | 0.00 |
| from * company | 44 | 1029.74 | 23.40 | 60.01 | 0.00 |
| from * hour | 1056 | 1793.35 | 1.70 | 4.35 | 0.00 |
| from * week | 44 | 426.40 | 9.69 | 24.85 | 0.00 |
| company * hour | 24 | 29.32 | 1.22 | 3.13 | 0.00 |
| company * week | 1 | 13.73 | 13.73 | 35.20 | 0.00 |
| hour * week | 24 | 43.20 | 1.80 | 4.62 | 0.00 |
| | | | | | |
| to * from * company | 132 | 487.21 | 3.69 | 9.46 | 0.00 |
| to * from * hour | 3168 | 1326.40 | 0.42 | 1.07 | 0.02 |
| to * from * week | 132 | 162.25 | 1.23 | 3.15 | 0.00 |
| to * company * hour | 72 | 38.60 | 0.54 | 1.37 | 0.02 |
| to * company * week | 3 | 6.54 | 2.18 | 5.59 | 0.00 |
| to * hour * week | 72 | 25.91 | 0.36 | 0.92 | 0.66 |
| from * company * hour | 1056 | 745.65 | 0.71 | 1.81 | 0.00 |
| from * company * week | 44 | 139.37 | 3.17 | 8.12 | 0.00 |
| from * hour * week | 1056 | 782.30 | 0.74 | 1.90 | 0.00 |
| company * hour * week | 24 | 24.51 | 1.02 | 2.62 | 0.00 |
| | | | | | |
| to * from * company * hour | 3168 | 1339.13 | 0.42 | 1.08 | 0.01 |
| to * from * company * week | 132 | 117.49 | 0.89 | 2.28 | 0.00 |
| to * from * hour * week | 3168 | 1308.72 | 0.41 | 1.06 | 0.05 |
| to * company * hour * week | 72 | 31.62 | 0.44 | 1.13 | 0.22 |
| from * company * hour * week | 1056 | 528.34 | 0.50 | 1.28 | 0.00 |
| | | | | | |
| to * from * company * hour * week | 3168 | 1235.54 | 0.39 | | |

Fig. 4.  *Classical ANOVA table for a $4 \times 45 \times 2 \times 25 \times 2$ factorial data structure. The data are logarithms of connect times for messages on the World Wide Web.*

"to" and "from" locations and the time ordering of the hours. As is usual, ANOVA is a tool for data exploration—for learning about which factors are important in predicting the variation in the data—which can be used to construct useful models or design future data collection. The linear model is a standard approach to analyzing factorial data; in this context, we see that the multilevel ANOVA display, which focuses on variance components, conveys more relevant information than does the classical ANOVA, which focuses on null hypothesis testing.

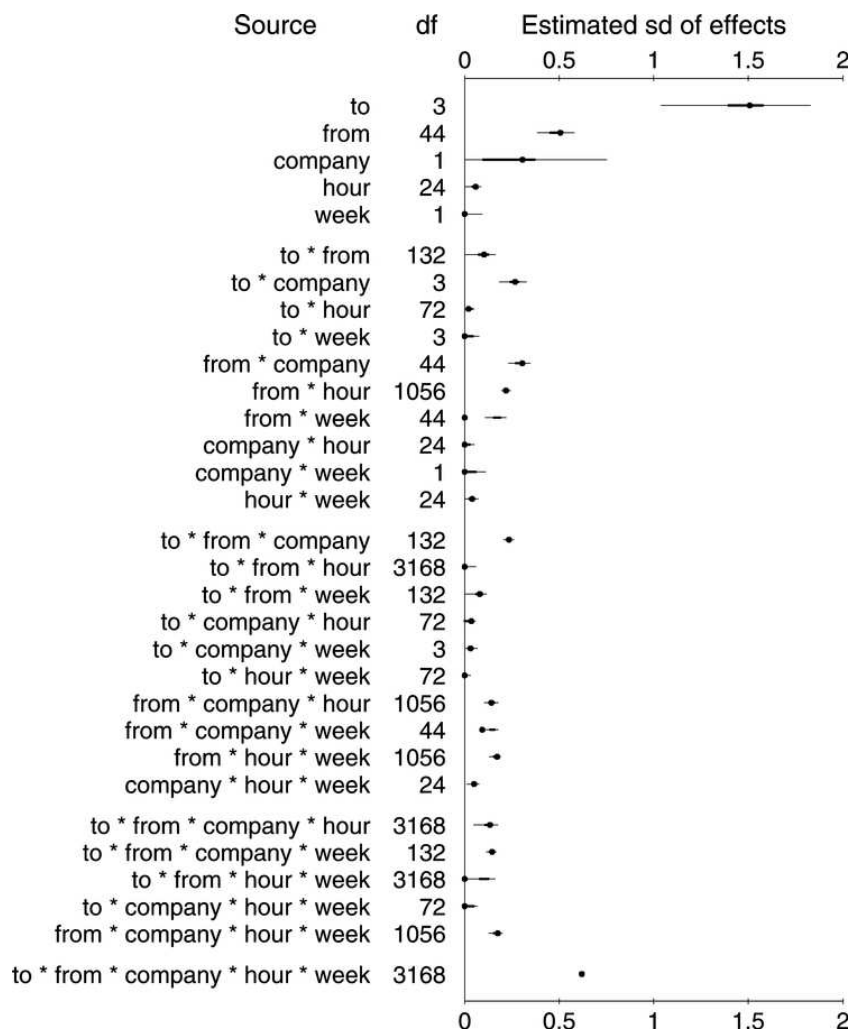

Fig. 5. *ANOVA display for the World Wide Web data (cf. to the classical ANOVA in Figure 4). The bars indicate 50% and 95% intervals for the finite-population standard deviations $s_m$, computed using simulation based on the classical variance component estimates. Compared to the classical ANOVA in Figure 4, this display makes apparent the magnitudes and uncertainties of the different components of variation. Since the data are on the logarithmic scale, the standard deviation parameters can be interpreted directly. For example, $s_m = 0.20$ corresponds to a coefficient of variation of $\exp(0.2) - 1 \approx 0.2$ on the original scale, and so the unlogged coefficients $\exp(\beta_j^{(m)})$ in this batch correspond to multiplicative increases or decreases in the range of 20%.*

Another direction to consider is the generalization of the model to new situations. Figure 5 displays uncertainty intervals for the finite-population standard deviations so as to be comparable to classical ANOVA. This makes sense when comparing the two companies and 25 hours, but the "to" sites, the "from" sites and the weeks are sampled from a larger population, and for these generalizations, the superpopulation variances would be relevant.

### 7.2. A multilevel logistic regression model with interactions: political opinions.
Dozens of national opinion polls are conducted by media organizations before every election, and it is desirable to estimate opinions at the levels of individual states as well as for the entire country. These polls are generally based on national random-digit dialing with corrections for nonresponse based on demographic factors such as sex, ethnicity, age and education [see Voss, Gelman and King (1995)]. We estimated state-level opinions from these polls, while simultaneously correcting for nonresponse, in two steps. For any survey response of interest:

1. We fit a regression model for the individual response given demographics and state. This model thus estimates an average response $\theta_j$ for each cross-classification $j$ of demographics and state. In our example, we have sex (male/female), ethnicity (black/nonblack), age (four categories), education (four categories) and 50 states; thus 3200 categories.
2. From the Census, we get the adult population $N_j$ for each category $j$. The estimated average response in any state $s$ is then $\theta_s = \sum_{j \in s} N_j \theta_j / \sum_{j \in s} N_j$, with each summation over the 64 demographic categories in the state.

We need a large number of categories because (a) we are interested in separating out the responses by state, and (b) nonresponse adjustments force us to include the demographics. As a result, any given survey will have few or no data in many categories. This is not a problem, however, if a multilevel model is fit, as is done automatically in our ANOVA procedure: each factor or set of interactions in the model, corresponding to a row in the ANOVA table, is automatically given a variance component.

As described by Gelman and Little (1997) and Bafumi, Gelman and Park (2002), this inferential procedure works well and outperforms standard survey estimates when estimating state-level outcomes. For this paper, we choose a single outcome—the probability that a respondent prefers the Republican candidate for President—as estimated by a logistic regression model from a set of seven CBS News polls conducted during the week before the 1988 Presidential election. We focus here on the first stage of the estimation procedure—the inference for the logistic regression model—and use our ANOVA tools to display the relative importance of each factor in the model.

We label the survey responses $y_i$ as 1 for supporters of the Republican candidate and 0 for supporters of the Democrat (with undecideds excluded) and model them as independent, with $\Pr(y_i = 1) = \text{logit}^{-1}((X\beta)_i)$. The design matrix $X$ is all 0's and 1's with indicators for the demographic variables used by CBS in the survey weighting: sex, ethnicity, age, education and the interactions of sex × ethnicity and age × education. We also include in $X$ indicators for the 50 states and for the four regions of the country (northeast, midwest, south and west). Since the states are nested within regions (which is implied by the design matrix of the regression), no main effects for states are needed. As in our general approach for linear models, we give each batch of regression coefficients an independent normal distribution centered at zero and with standard deviation estimated hierarchically given a uniform prior density.

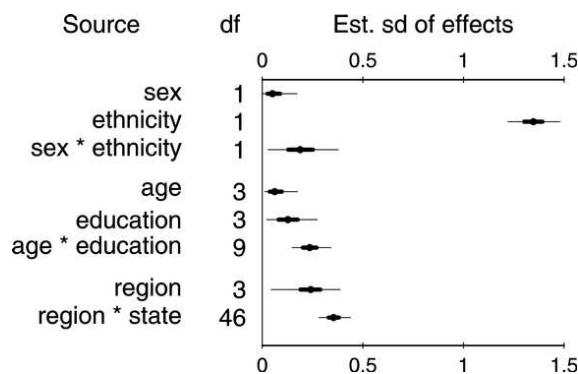

Fig. 6.   *ANOVA display for the logistic regression model of the probability that a survey respondent prefers the Republican candidate for the 1988 U.S. Presidential election, based on data from seven CBS News polls. Point estimates and error bars show posterior medians, 50% intervals and 95% intervals of the finite-population standard deviations $s_m$, computed using Bayesian posterior simulation. The demographic factors are those used by CBS to perform their nonresponse adjustments, and states and regions are included because we were interested in estimating average opinions by state. The large effects for ethnicity and the general political interest in states suggest that it might make sense to include interactions; see Figure 7.*

We fit the model using the Bayesian software Bugs [Spiegelhalter, Thomas, Best and Lunn (2002)], linked to R [R Project (2000) and Gelman (2003)] where we computed the finite-sample standard deviations and plotted the results. Figure 6 displays the ANOVA table, which shows that ethnicity is by far the most important demographic factor, with state also explaining quite a bit of variation.

The natural next step is to consider interactions among the most important effects, as shown in Figure 7. The `ethnicity * state * region` interactions are surprisingly large: the differences between African-Americans and others vary dramatically by state. As with the previous example, ANOVA is a useful tool in understanding the importance of different components of a hierarchical model.

**8. Discussion.** In summary, we have found hierarchical modeling to be a key step in allowing ANOVA to be performed reliably and automatically. Conversely, the ideas of ANOVA are extremely powerful in modeling complex data of the sort that we increasingly handle in statistics—hence the title of this paper. We conclude by reviewing these points and noting some areas for further work.

8.1. *The importance of hierarchical modeling in formulating and computing ANOVA.* Analysis of variance is fundamentally about multilevel modeling: each row in the ANOVA table corresponds to a different batch of parameters, along with inference about the standard deviation of the parameters in this batch. A crucial difficulty in classical ANOVA and, more generally, in classical linear modeling, is identifying the correct variance components to use in computing standard errors and testing hypotheses. The hierarchical data structures in Section 2.2 illustrate the limitations of performing ANOVA using classical regression.

However, as we discuss in this paper, assigning probability distributions for all variance components automatically gives the correct comparisons and standard errors. Just as a design matrix corresponds to a particular linear model, an ANOVA table corresponds to a particular multilevel batching of random effects. It should thus be possible to fit any ANOVA automatically without having to figure out the appropriate error variances, even for notoriously difficult designs such as split-plots (recall Figure 1).

8.2. *Estimation and hypothesis testing for variance components.* This paper has identified ANOVA with estimation in variance components models. As discussed in Section 3.5, uncertainties can be much lower for finite-population variances $s_m^2$ than for superpopulation variances $\sigma_m^2$, and it is through finite-population variances that we connect to classical ANOVA, in which it is possible to draw useful inferences for even small batches (as in our split-plot Latin square example).

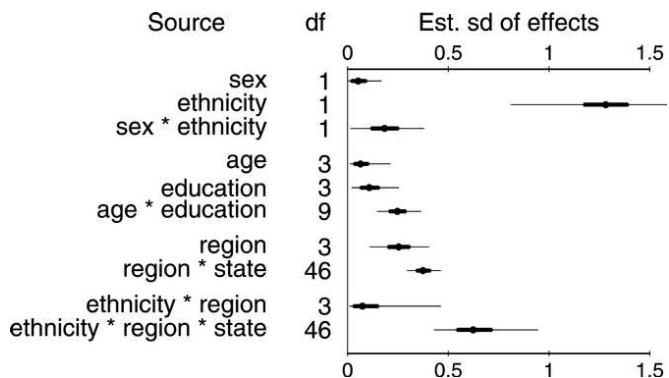

FIG. 7. *ANOVA display for the logistic regression model for vote preferences, adding interactions of ethnicity with region and state. Compare to Figure 6.*

Hypothesis testing is in general a more difficult problem than estimation because many different possible hypotheses can be considered. In some relatively simple balanced designs, the hypotheses can be tested independently; for example, the split-plot Latin square allows independent testing of row, column and treatment effects at the between- and within-plot levels. More generally, however, the test of the hypothesis that some $\sigma_m = 0$ will depend on the assumptions made about the variance components lower in the table. For example, in the factorial analysis of the Internet data in Section 7.1, a test of the `to * from` interaction will depend on the estimated variances for all the higher-level lower interactions including `to * from`, and it would be inappropriate to consider only the full five-way interaction as an "error term" for this test (since, as Figures 4 and 5 show, many of the intermediate outcomes are both statistically significant and reasonably large). Khuri, Mathew and Sinha (1998) discuss some of the options in testing for variance components, and from a classical perspective these options proliferate for unbalanced designs and highly structured models.

From a Bayesian perspective, the corresponding step is to model the variance parameters $\sigma_m$. Testing for null hypotheses of zero variance components corresponds to hierarchical prior distributions for the variance components that have a potential for nonnegligible mass near zero, as has been discussed in the Bayesian literature on shrinkage and model selection [e.g., Gelman (1992), George and McCulloch (1993) and Chipman, George and McCulloch (2001)]. In the ANOVA context such a model is potentially more difficult to set up since it should ideally reflect the structure of the variance components (e.g., if two sets of main effects are large, then one might expect their interaction to be potentially large).

8.3. *More general models.* Our model (7) for the linear parameters corresponds to the default inferences in ANOVA, based on computations of variances and exchangeable coefficients within each batch. This model can be expanded in various ways. Most simply, the distributions for the effects in each batch can be generalized beyond normality (e.g., using $t$ or mixture distributions), and the variance parameters can themselves be modeled hierarchically, as discussed immediately above.

Another generalization is to nonexchangeable models. A common way that nonexchangeable regression coefficients arise in hierarchical models is through group-level regressions. For example, the five rows, columns and possibly treatments in the Latin square are ordered, and systematic patterns there could be modeled, at the very least, using regression coefficients for linear trends. In the election survey example, one can add state-level predictors such as previous Presidential election results. After subtracting batch-level regression predictors, the additive effects for the factor levels in each batch could be modeled as exchangeable. This corresponds to analysis of covariance or contrast analysis in classical ANOVA. Our basic model (6) sets up a regression at the level of the data, but regressions on the hierarchical coefficients (i.e., contrasts) can have a different substantive interpretation as interblock or contextual effects [see Kreft and de Leeuw (1998) and Snijders and Bosker (1999)]. In either case, including contrasts adds another twist in that defining a superpopulation for predictive purposes now requires specifying a distribution over the contrast variable (e.g., in the Latin square example, if the rows are labeled as $-2, -1, 0, 1, 2$, then a reasonable superpopulation might be a uniform distribution on the range $[-2.5, 2.5]$).

More complex structures, such as time-series and spatial models [see Ripley (1981) and Besag and Higdon (1999)], or negative intraclass correlations, cannot be additively decomposed in a natural way into exchangeable components. One particularly interesting class of generalizations of classical ANOVA involves the nonadditive structures of interactions. For example, in the Internet example in Section 7.1 the coefficients in any batch of two-way or higher-level interactions have a natural gridded structure that is potentially more complex than the pure exchangeability of additive components [see Aldous (1981)].

8.4. *The importance of the ANOVA idea in statistical modeling and inference.* ANOVA is more important than ever because it represents a key idea in statistical modeling of complex data structures—the grouping of predictor variables and their coefficients into batches. Hierarchical modeling, along with the structuring of input variables, allows the modeler easily to include hundreds of predictors in a regression model (as with the examples in Section 7), as has been noted by proponents of multilevel modeling [e.g., Goldstein (1995), Kreft and de Leeuw (1998) and Snijders and Bosker (1999)]. ANOVA allows us to understand these models in a way that we cannot by simply looking at regression coefficients, by generalizing classical variance components estimates [e.g., Cochran and Cox (1957) and Searle, Casella and McCulloch (1992)]. The ideas of the analysis of variance also help us to include finite-population and superpopulation inferences in a single fitted model, hence unifying fixed and random effects. A future research challenge is to generalize our inferences and displays to include multivariate models of coefficients (e.g., with random slopes and random intercepts, which will jointly have a covariance matrix as well as individual variances).

**Acknowledgments.** We thank Hal Stern for help with the linear model formulation; John Nelder, Donald Rubin, Iven Van Mechelen and the editors and referees for helpful comments; and Alan Edelman for the data used in Section 7.1.

Department of Statistics
Columbia University
New York, New York 10027
USA
e-mail: gelman@stat.columbia.edu